\documentclass[onethmnum,oneeqnum,final]{siamonline190516}

\usepackage{amsmath}
\usepackage{amsfonts}
\usepackage{fullpage}
\usepackage{cleveref}

\newcommand{\CC}{\ensuremath{\mathbb{C}}}
\newcommand{\RR}{\ensuremath{\mathbb{R}}}
\newcommand{\QQ}{\ensuremath{\mathbb{Q}}}
\newcommand{\ZZ}{\ensuremath{\mathbb{Z}}}
\newcommand{\V}{\ensuremath{\mathcal{V}}}

\newcommand{\boldalpha}{\ensuremath{\boldsymbol{\alpha}}}
\newcommand{\hatalpha}{\ensuremath{\hat{\boldsymbol{\alpha}}}}
\newcommand{\rowv}[1]{\ensuremath{\mathbf{#1}}}
\newcommand{\colv}[1]{\ensuremath{\boldsymbol{#1}}}
\newcommand{\pt}[1]{\ensuremath{\boldsymbol{#1}}}
\newcommand{\inner}[2]{\ensuremath{ \left\langle #1 \,,\, #2 \right\rangle }}

\DeclareMathOperator{\nvol}{V}
\DeclareMathOperator{\mvol}{MV}
\DeclareMathOperator{\conv}{conv}
\DeclareMathOperator{\init}{init}
\DeclareMathOperator{\supp}{supp}

\title{Volume of convex polytopes equals mixed volume of simplices}

\author{Tianran Chen\thanks{Department of Mathematics, Auburn University at Montgomery, Montgomery Alabama USA}
    \email{ti@nranchen.org}
}

\begin{document}

\maketitle

\begin{abstract}
    This note provides a simple proof for the equality between 
    the normalized volume of a convex polytope with $m$ vertices
    and the mixed volume of $m$ simplices
    and thus shows the seemingly restrictive problem of computing mixed volume
    of simplices is still at least as hard as computing volumes of convex polytopes.
\end{abstract}



\section{Introduction}

Let $\nvol$ and $\mvol$ denote the normalized and mixed volume functions
for lattice polytopes respectively. 
It is well known that for an $n$-simplex spanned by
$\pt{p}_1, \dots, \pt{p}_{n+1} \in \RR^n$,
\[
    \nvol(\conv\{\, \pt{p}_1, \dots, \pt{p}_{n+1} \}) =
    \left|
    \det
    \begin{bmatrix}
        1 & \cdots & 1 \\
        \pt{p}_1 & \cdots & \pt{p}_{n+1}
    \end{bmatrix}
    \right|
    =
    \mvol(
        \{ \hat{\pt{p}}_1,     \pt{e}_{n+1} \},
        \dots,
        \{ \hat{\pt{p}}_{n+1}, \pt{e}_{n+1} \}
    ),
\]
where $\pt{p}_i$'s are expressed as column vectors,
$\hat{\pt{p}}_i$ is the embedding of $\pt{p}_i$ into
$\RR^n \times \{ \pt{0} \} \subset \RR^m$,
and $\pt{e}_{n+1} \in \RR^{n+1}$ is the $(n+1)$-th unit vector.
This equation tells us that
the normalized volume of an $n$-simplex in $\RR^n$
can be turned into the mixed volume of $n+1$ line segments.
Using root counting results from algebraic geometry,
this short notes provides a simple proof for the generalization of this observation
to cases of convex hull of $m$ points $\pt{p}_1,\dots,\pt{p}_m$ in $\RR^n$ for $m > n$:
\[
    \nvol ( \conv\{\, \pt{p}_1, \dots, \pt{p}_m \}) =
    \mvol ( \Delta(\pt{p}_1), \dots, \Delta(\pt{p}_m) )
\]
where each $\Delta(\pt{p}_i)$ is a simplex of dimension $m-n$.
From the view point of complexity,
this equation shows that the seemingly specialized problem of computing
the mixed volume of simplices is still at least as hard as
the problem of computing the volume of a convex polytope.

\section{Notations and preliminaries}

$\CC,\RR,\QQ,\ZZ$ are the sets of
complex numbers, reals, rationals, and integers respectively.
The set of $m \times n$ matrices with entries in a ring $R$ is denoted
$M_{m \times n}(R)$ 

For a set $S \subseteq \RR^n$, $\conv(S)$ is the convex hull of $S$,
i.e., the smallest convex set containing $S$.
A \emph{convex polytope} in $\RR^n$ is a subset $P$ such that 
$P = \conv\{\, \pt{p}_1, \dots, \pt{p}_m\}$
for some $\pt{p}_1, \dots, \pt{p}_m \in \RR^n$.
Its dimension $\dim(P)$ is the dimension of the smallest affine space containing it,
and its \emph{normalized volume} $V(P)$ is $n! \operatorname{Vol}_n(P)$
where $\operatorname{Vol}_n$ is the Euclidean volume form for $\RR^n$.
Note that this definition is not restricted to lattice polytopes
and implies the convention that $\nvol(P) = 0$ if $\dim(P) < n$.

For two sets $A,B \subset \RR^n$, their \emph{Minkowski sum} is
$A+B = \{ \pt{a} + \pt{b} \mid \pt{a} \in A , \pt{b} \in B \}$.
For convex polytopes $P_1,\dots,P_n \subset \RR^n$
the volume of the Minkowski sum $\lambda_1 P_1 + \cdots + \lambda_n P_n$
under the scaling by positive factors $\lambda_1,\dots,\lambda_n$,
is a homogeneous polynomial in $\lambda_1,\dots,\lambda_n$~\cite{Minkowski1911Theorie}.
The coefficient of the mixed term $\lambda_1 \cdots \lambda_n$ is known as 
the \emph{mixed volume} of these polytopes, denoted $\mvol(P_1,\dots,P_n)$.

For $\rowv{x} = \begin{bmatrix} \, x_1 & \dots & x_n \end{bmatrix}$
and $\colv{a} = \begin{bmatrix} \, a_1 & \dots & a_n \end{bmatrix}^\top \in \ZZ^n$,
$\rowv{x}^{\colv{a}} = x_1^{a_1} \, \cdots \, x_n^{a_n}$.
Similarly, for $A \in M_{n \times m}(\ZZ)$
with columns $\colv{a}_1,\dots,\colv{a}_m \in \ZZ^n$,
$ \rowv{x}^A = \begin{bmatrix} \rowv{x}^{\,\colv{a}_1} & \dots & \rowv{x}^{\colv{a}_m} \end{bmatrix}$,
which defines functions over $(\CC^*)^n = (\CC \setminus \{ 0 \})^n$.
A \emph{Laurent polynomial} is an expression of the form
$f(\rowv{x}) = \sum_{\colv{a} \in S} c_{\colv{a}} \rowv{x}^{\colv{a}}$
where each $c_{\colv{a}} \in \CC^*$,
$S \subset \ZZ^n$ is the the \emph{support} of $f$, denoted $\supp(f)$,
and $\conv(S)$ is its \emph{Newton polytope}.
For a system of Laurent polynomials $F = (f_1,\dots,f_n)$ in 
$\rowv{x} = \begin{bmatrix} \, x_1 & \dots & x_n \end{bmatrix}$,
its set of common isolated zeros in $(\CC^*)^n$ is denoted $\V_0^*(F)$,
and the following theorems provide us bound on $|\V_0^*(F)|$.

\begin{theorem}[Kushnirenko \cite{Kushnirenko1975Newton}]
    \label{thm: Kushnirenko}
    If $S := \supp(f_i)$ for all $i$ are identical,
    $|\V_0^*(F)| \le \nvol (\conv(S))$.
\end{theorem}

\begin{theorem}[Bernshtein's 1st Theorem \cite{Bernshtein1975Number}]
    \label{thm: Bernshtein1}
    $|\V_0^*(F)| \le \mvol (\conv(\supp(f_1)),\dots,\conv(\supp(f_n)))$.
\end{theorem}

Moreover, for generic choices of the coefficients,
both bounds are exact.
In general, 
the exactness of this bound can be tested via
``initial systems''.
For a Laurent polynomial 
$f(\rowv{x}) = \sum_{\colv{a} \in S} c_{\colv{a}} \rowv{x}^{\colv{a}}$
and $\boldalpha \in \RR^n$,
$\init_{\boldalpha}(f)(\rowv{x}) = 
\sum_{\colv{a} \in (S)_{\boldalpha}} c_{ \colv{a}} \rowv{x}^{\colv{a}}$,
where $(S)_{\boldalpha}$ is the subset of $S$ on which the linear functional
$\inner{ \boldalpha }{ \cdot }$ is minimized.
For a Laurent polynomial system $F = (f_1,\dots,f_n)$ in $x_1,\dots,x_n$,
its \emph{initial system} $\init_{\boldalpha} (F)$ with respect a vector
$\boldalpha \in \RR^n$, is 
$(\init_{\boldalpha}(f_1),\dots,\init_{\boldalpha}(f_1))$.

\begin{theorem}[Bernshtein's 2nd Theorem \cite{Bernshtein1975Number}]
    \label{thm: Bernshtein2}
    If for all nonzero vectors $\boldalpha \in \RR^n$,
    the initial system $\init_{\boldalpha}(F)$ has no zero in $(\CC^*)^n$,
    then all zeros of $F$ in $(\CC^*)^n$ are isolated,
    and the total number, counting multiplicity, is 
    $\mvol(\conv(S_1),\dots,\conv(S_n))$.
    The converse is also true.
\end{theorem}

We avoid defining ``multiplicity''
as only zeros of multiplicity 1 (nonsingular zeros) are involved.


\section{The main result}

Using the root counting theorems listed above,
collectively known as the theory of BKK bound \cite{Bernshtein1975Number,Khovanskii1978Newton,Kushnirenko1975Newton}
we now state and prove the main result.

\begin{theorem}\label{thm: main}
    For $m$ distinct points $\{ \pt{p}_1, \dots, \pt{p}_m \} \subset \RR^n$,
    with $m > n > 0$,
    \[
        \nvol ( \conv\{\, \pt{p}_1, \dots, \pt{p}_m \}) =
        \mvol ( \Delta(\pt{p}_1), \dots, \Delta(\pt{p}_m) )
    \]
    where
    \begin{align*}
        \Delta(\pt{p}) &= \conv\{ \, \hat{\pt{p}}, \, \pt{e}_{n+1}, \dots, \pt{e}_m \}.
        &&\text{and} &
        \hat{\pt{p}} = 
        \begin{bmatrix}
            \pt{p} \\
            \pt{0}_{m-n}
        \end{bmatrix}
        \in \RR^m.
    \end{align*}
\end{theorem}

Here, $\Delta(\pt{p})$ is a simplex of dimension $m-n$.
So the equation above states that the normalized volume
of a finite set of the convex hull of $m$ points can always be expressed as
the mixed volume of $m$ simplices in $\RR^m$.
Also recall convention that
$\nvol (X) = 0$ if and only if $X$ is not full-dimensional.

\begin{proof}
    In the trivial case where $\dim(\conv \{\, \colv{p}_1,\dots,\colv{p}_m \}) < n$,
    we adopted the convention that this is equivalent to
    $\nvol(\conv \{\, \colv{p}_1,\dots,\colv{p}_m \}) = 0$.
    The simplices $\Delta(\pt{p}_i)$ for $i=1,\dots,m$ are all contained in
    a common hyperplane in $\RR^m$, and thus
    $\mvol ( \Delta(\pt{p}_1), \dots, \Delta(\pt{p}_m) )$ is also zero.
    Conversely, if $\mvol ( \Delta(\pt{p}_1), \dots, \Delta(\pt{p}_m) ) = 0$,
    then there is no $m$-tuple of line segments $(L_1,\dots,L_m)$ with
    $L_i \subset \Delta(\pt{p}_i)$ for $i=1,\dots,m$ that represent linearly independent vectors.
    In that case, $\{\, \pt{p}_1, \dots, \pt{p}_m \}$ must be affinely dependent,
    and hence $\nvol(\conv \{\, \colv{p}_1,\dots,\colv{p}_m \}) = 0$.

    For cases where
    $\conv \{\, \pt{p}_1,\dots,\pt{p}_m \}$ is full-dimensional,
    we first assume $\{\, \pt{p}_1, \dots, \pt{p}_m \} \subset \ZZ^n$.
    Let
    \[
        P = 
        \begin{bmatrix} 
            \; \pt{p}_1 & \cdots & \pt{p}_m \;
        \end{bmatrix}
        \in M_{n \times m}(\ZZ),
    \]
    then for $\rowv{x} = \begin{bmatrix} \, x_1 & \cdots & x_n \, \end{bmatrix}$
    with $x_1,\dots,x_n \in \CC^*$,
    $\rowv{x}^P$ is the system of $m$ Laurent monomials
    $\begin{bmatrix} \,\rowv{x}^{\pt{p}_1} & \cdots & \rowv{x}^{\pt{p}_m} \,\end{bmatrix}$
    expressed as a row vector.
    For a generic choice of an $n \times m$ complex matrix $A$,
    we consider the square system of $n$ Laurent polynomials
    \[
        F(\rowv{x}) = A \, (\rowv{x}^P)^\top
    \]
    in $\rowv{x}$.
    By the genericity assumption,
    the zero set of $F$ in $(\CC^*)^n$ consists of nonsingular isolated points,
    and, by Kushnirenko's Theorem (\Cref{thm: Kushnirenko}),
    $|\V_0^*(F)| = \nvol( \conv \{\, \pt{p}_1, \dots, \pt{p}_m \})$.

    Since $m > n$, and the choice of $A \in M_{n \times m}(\CC)$ is generic,
    the dimension of the null space of $A$ is exactly $d = m - n > 0$.
    Let $\{ \colv{v}_1,\dots,\colv{v}_d \}$ be a basis of the null space of $A$
    and $K = \begin{bmatrix} \colv{v}_1 & \cdots & \colv{v}_d \end{bmatrix} \in M_{m \times d}(\CC)$,
    then, for $\rowv{x} \in (\CC^*)^n$,
    \begin{align*}
        F(\rowv{x}) &= A \, (\rowv{x}^P)^\top = \colv{0}
        &&\text{if and only if} &
        (\rowv{x}^P)^\top &= K (\rowv{y})^\top 
    \end{align*}
    for some $\rowv{y} = \begin{bmatrix} \,y_1 & \cdots & y_d\,\end{bmatrix} \in \CC^d$.
    Moreover, since $\V_0^*(F)$ is finite, without loss of generality, we can choose the basis
    $\colv{v}_1,\dots,\colv{v}_d$ so that $\rowv{y} \in (\CC^*)^d$
    for every corresponding $\rowv{x} \in \V_0^*(F)$.
    Define $G = (g_1,\dots,g_m) : (\CC^*)^n \times (\CC^*)^d \to \CC^m$ given by
    \[
        G(\rowv{x},\rowv{y}) = (\rowv{x}^P)^\top - K (\rowv{y})^\top.
    \]
    Then $G$ is a system of nonzero Laurent polynomials in $x_1,\dots,x_n,y_1,\dots,y_d$,
    and there is a bijection between $\V_0^*(F)$ and $\V_0^*(G)$.
    Moreover, under the genericity assumption, this map preserves multiplicity,
    and therefore we simply need to show 
    $|\V_0^*(G)| = \mvol ( \Delta(\pt{p}_1), \dots, \Delta(\pt{p}_m) )$.

    The Newton polytopes of $G$ are
    $( \Delta(\pt{p}_1), \dots, \Delta(\pt{p}_m) )$.
    Therefore, by Bernshtein's First Theorem,
    \[ |\V_0^*(G)| \le \mvol ( \Delta(\pt{p}_1), \dots, \Delta(\pt{p}_m) ). \]
    To establish equality, we shall show the system $G$ satisfies the conditions
    in \Cref{thm: Bernshtein2},
    even though there are algebraic relations among the coefficients.
    That is, we aim to show for any nonzero vector
    $\hat{\boldalpha} \in \RR^m$,
    the initial system $\init_{\hat{\boldalpha}} G$ has no $\CC^*$-zeros.
    Fix a nonzero $\hat{\boldalpha} = (\alpha_1,\dots,\alpha_m) \in \RR^m$.
    We define $\boldalpha = (\alpha_1,\dots,\alpha_n)$,
    $\ell = \min \{ \alpha_{n+j} \mid j = 1,\dots,d \}$ and
    $L = \{ j \in \{ 1,\dots, d \} \mid \alpha_{n+j} = \ell \}$,
    then 
    \begin{align*}
        \init_{\hatalpha}(g_i)(\rowv{x},\rowv{y}) &= 
        \begin{cases}
            \rowv{x}^{\colv{p}_i} 
            &\text{if } \inner{ \boldalpha }{ \colv{p}_i } < \ell \\
            \rowv{x}^{\colv{p}_i}
            - \sum_{j \in L} v_{ij} y_{j}
            &\text{if } \inner{ \boldalpha }{ \colv{p}_i } = \ell \\
            \phantom{\rowv{x}^{\colv{p}_i}}
            - \sum_{j \in L} v_{ij} y_{j}
            &\text{if } \inner{ \boldalpha }{ \colv{p}_i } > \ell \\
        \end{cases}
        &&\text{for each } i=1,\dots,m.
    \end{align*}

    If $\inner{ \boldalpha }{ \colv{p}_i } < \ell$,
    i.e., the first case appears, for any $i \in \{1,\dots,m\}$,
    then $\init_{\hatalpha}(G)$ has no $\CC^*$-zero.
    It is therefore sufficient to only consider cases where
    $\inner{ \boldalpha }{ \colv{p}_i } \ge \ell$ for all $i$.
    Furthermore, since we assumed $\rowv{y} \in (\CC^*)^d$ and is unique
    for every corresponding $\rowv{x} \in \V_0^*(F)$,
    it is sufficient to only consider cases where
    $\alpha_{n+1},\dots,\alpha_{m} \ge 0$.

    Let
    $B = \{ i \in \{1,\dots,m\} \mid \inner{ \boldalpha }{ \colv{p}_i } = \ell \}$.
    Since $\conv\{\, \colv{p}_1, \dots, \colv{p}_m \}$ is
    assumed to be full-dimensional, we must have $B \ne \{1,\dots,m\}$
    and hence $\boldalpha \ne \rowv{0}$,
    In this case,
    any $(\rowv{x},\rowv{y}) \in (\CC^*)^n \times (\CC^*)^d$
    satisfying $\init_{\hatalpha}(G)(\rowv{x},\rowv{y}) = \colv{0}$
    gives raise to a $\CC^*$-zero to
    $\init_{\boldalpha}(F)(\rowv{x}) = \colv{0}$,
    which contradict with the assumption that
    $|\V_0^*(F)| = \nvol(\conv\{\, \colv{p}_1, \dots, \colv{p}_m \})$,
    according to \Cref{thm: Bernshtein2}.
    Therefore, we can conclude that $\init_{\hatalpha}(G)(\rowv{x},\rowv{y})$
    has no $\CC^*$-zeros.

    We have shown that $\init_{\hatalpha}(G)$ has no $\CC^*$ zeros
    for any nonzero $\hatalpha \in \RR^n$.
    Therefore, by \Cref{thm: Bernshtein2},
    \[
        \nvol ( \conv\{\, \pt{p}_1, \dots, \pt{p}_m \}) =
        | \V^* (F) | =
        | \V^* (G) | = 
        \mvol ( \Delta(\pt{p}_1), \dots, \Delta(\pt{p}_m) )
    \]
    for a set $\{\, \colv{p}_1, \dots, \colv{p}_m \}$ of distinct points in $\ZZ^n$.

    Since both 
    $\nvol ( \conv\{\, \pt{p}_1, \dots, \pt{p}_m \})$ and
    $\mvol ( \Delta(\pt{p}_1), \dots, \Delta(\pt{p}_m) )$
    are homogeneous of degree $n$ with respect to a uniform scaling
    $\pt{p}_i \mapsto \lambda \pt{p}_i$,
    the equality extends to cases where
    $\{\, \colv{p}_1, \dots, \colv{p}_m \} \subset \QQ^n$.
    By continuity with respect to the Hausdorff metric,
    the result further extends to point sets in $\RR^n$.
\end{proof}

\section{Complexity implications}

It is well known that the mixed volume of $n$ line segments in $\RR^n$
is equivalent to the determinant of an $n \times n$ matrix
and hence can be computed in polynomial time.
Yet, M. Dyer, P. Gritzmann, and A. Hufnagel
showed that even the relatively simple problem of computing the mixed volume of $n$ ``boxes'',
each formed as the Minkowski sum of $n$ line segments of nonuniform length along the axes,
is $\#\mathbb{P}$-hard \cite[Theorem 3]{dyer_complexity_1998}.
It is reasonable to speculate that one source of the $\#\mathbb{P}$-hardness
is the geometric complexity of the ``boxes'',
since a $d$-dimensional box can be subdivided into $d!$ full-dimensional simplices.
This is not the case.
\Cref{thm: main} shows that even when the ``boxes'' are replaced by simplices,
the $\#\mathbb{P}$-hardness persists:
Even the much more restrictive problem of 
computing the mixed volume of $m$ simplices of dimension $d < m$
is at least as hard as the problem of computing the volume of the
convex hull of $m$ points in $\RR^{m-d}$,
which is known to be $\#\mathbb{P}$-hard.

\bibliographystyle{siamplain}
\bibliography{library}

\begin{thebibliography}{1}

\bibitem{Bernshtein1975Number}
{\sc D.~N. Bernshtein}, {\em {The number of roots of a system of equations}},
  Functional Analysis and its Applications, 9 (1975), pp.~183--185.

\bibitem{dyer_complexity_1998}
{\sc M.~Dyer, P.~Gritzmann, and A.~Hufnagel}, {\em {On the complexity of
  computing mixed volumes}}, SIAM Journal on Computing, 27 (1998),
  pp.~356--400, \url{https://doi.org/10.1137/S0097539794278384}.

\bibitem{Khovanskii1978Newton}
{\sc A.~G. Khovanskii}, {\em {Newton polyhedra and the genus of complete
  intersections}}, Functional Analysis and Its Applications, 12 (1978),
  pp.~38--46, \url{https://doi.org/10.1007/BF01077562},
  \url{http://dx.doi.org/10.1007/BF01077562}.

\bibitem{Kushnirenko1975Newton}
{\sc A.~G. Kushnirenko}, {\em {A Newton polyhedron and the number of solutions
  of a system of k equations in k unknowns}}, Usp. Math. Nauk, 30 (1975),
  pp.~266--267.

\bibitem{Minkowski1911Theorie}
{\sc H.~Minkowski}, {\em {Theorie der konvexen Korper, insbesondere Begrundung
  ihres Oberflachenbegriffs}}, Gesammelte Abhandlungen von Hermann Minkowski, 2
  (1911), pp.~131--229.

\end{thebibliography}

\end{document}